\documentclass[12pt, reqno]{amsart}
\usepackage{amsmath, amsthm, amscd, amsfonts, amssymb, graphicx, color}
\usepackage[bookmarksnumbered, colorlinks, plainpages]{hyperref}
\hypersetup{colorlinks=true,linkcolor=red, anchorcolor=green, citecolor=cyan, urlcolor=red, filecolor=magenta, pdftoolbar=true}
\usepackage{float}
\usepackage{multirow}
\usepackage{stackengine}
\usepackage{pdflscape}
\usepackage{tikz}
\usepackage{caption}
\usepackage{graphicx} 
\usepackage{hyperref}
\usetikzlibrary{decorations.pathreplacing}

\usepackage{setspace}  

\begin{document}
\title[]{Shape Gradient Based Non-Parametric Mumford-Shah Segmentation Without Level Sets}
\author[Shafeequdheen P]{Shafeequdheen P, Jyotiranjan Nayak, Vijayakrishna Rowthu}
\address{Department of Mathematics, SRM University-AP, Amaravati 522 502, India}
\email{shafeequdheen\_p@srmap.edu.in}
\begin{abstract}
A non-parametric, level-set-free method is proposed for detecting image boundaries using the shape gradient of the Mumford–Shah energy for segmentation. Minimizing the variance in pixel intensities inside and outside 
a boundary set of points is the primary pursuit. The boundary set as a polygon of points rather than a parametric form or a level set, evolves under the guidance of a shape gradient of the Mumford–Shah piece-wise constant segments model. Iteratively updating through the gradient descent method. The proposed method has been tested on various images, demonstrating its effectiveness in capturing intricate and narrow boundaries texture images.
\end{abstract}
\maketitle

\section{Introduction}
Image segmentation is a fundamental task in computer vision and image processing with applications like spanning object detection, medical imaging, and pattern recognition. The accurate segmentation enables meaningful interpretation and analysis of images by partitioning them into distinct regions based on specific features such as intensity, texture, or edges. Among the diverse segmentation techniques, variational models have garnered significant attention due to their mathematical rigor and robustness in handling noise, complex textures, and irregular structures~\cite{minaee2021image,yanowitz1989new,cheng2001color,zhang1996survey}. A well-known and widely used variational approach is the Mumford–Shah model~\cite{vitti2012mumford}, which formulates segmentation as an optimization problem, striving to divide an image into piecewise smooth regions while preserving essential edge information.

Building upon the variational segmentation framework, the work of E. Déreuve et al.~\cite{debreuve2007using,aubert2003image,barlaud2003shape} sheds light on the challenges posed by discretization in active contour methods, particularly when employing $L^2$-gradient descent for contour evolution. Their research demonstrated that the negativity condition of the shape derivative, a crucial requirement for stable evolution, is not always preserved in discrete settings, potentially leading to convergence issues. To mitigate this, they proposed a constrained approach wherein the velocity field is selected as a linear combination of predefined velocities, ensuring a more stable and flexible boundary evolution.

A general variational approach to segmentation involves defining an energy functional that captures local or global image characteristics, with segmentation achieved by minimizing this energy. The shape derivative of the energy functional~\cite{sokolowski1992introduction,cai2013two,chicco2017shape,delfour2011shapes} plays a pivotal role in this process, as it determines the velocity field that deforms the domain boundary to reduce the energy. This iterative deformation process is commonly referred to as the active contour method, which has been extensively used for boundary extraction and object delineation.

For 2D image segmentation, the Chan–Vese model~\cite{vese2002multiphase,huang2021chan,ruiying2022method}, a level-set-based active contour method, provides a simplification of the Mumford–Shah formulation by evolving a contour based on intensity differences inside and outside a given region. This process eliminates the reliance on explicit edge detection, level-set methods introduce additional computational overhead due to re-initialization requirements and numerical stability constraints~\cite{lie2006variant}. These factors can make level-set-based approaches computationally demanding, particularly for high-resolution images. Moreover, the contours generated from level sets are well-known for their significant changes in topology. However, this may not always be desirable for certain segmentation tasks, especially when preserving delicate boundaries without altering the topology is essential.

In this work, we propose a non-parametric level-set-free approach that leverages the shape derivative of the energy functional for 2D image segmentation by enhancing it while maintaining theoretical consistency, particularly in cases where prior knowledge about the object's structure can be integrated into the process. Instead of relying on level-set methods, our approach minimizes the variance of image intensity using a shape optimization framework, where the initial polygon of points  keeps evolving until convergence. The proposed method is validated on various images, including binary, RGB, and LAB color space images, demonstrating its effectiveness in capturing intricate structures while maintaining computational efficiency. By circumventing the challenges associated with level-set methods, our approach offers a robust alternative for high-resolution image segmentation tasks without altering the topology of the initial boundary.

The rest of this paper is organized as follows. Section \ref{Mathematical Model} introduces the mathematical foundation of the Mumford–Shah model. In Section \ref{Expression of the Shape Derivative}, we delve into shape derivatives, discussing examples of shape functionals and their derivatives, including the shape derivative of the Mumford–Shah model. Section \ref{Experimental validation} focuses on experimental validation, demonstrating the effectiveness of our proposed approach. Section \ref{Result analysis} provides a detailed analysis of the results.
\section{Mathematical Model}\label{Mathematical Model}
Let \(D \subseteq \mathbb{R}^2\) be a closed domain, and let \( \Omega \subset D \) be a simple closed region within \( D \) with boundary \( \Gamma \). The task is to partition the domain \( D \) into two regions: \( \Omega \) (the segmented region or the foreground) and its complement \( \Omega^c \) (the background). To achieve the segmentation, define an energy functional \( E(\Omega) \) as follows:
\begin{equation}\label{Mumford-shah model}
E(\Omega) =\frac{1}{|\Omega|} \int_{\Omega} (f(x) - \mu(\Omega))^2 \, dx+\frac{1}{|\Omega^c|} \int_{\Omega^c} (f(x) - \mu(\Omega^c))^2 \, dx+\eta\ \int_{\Gamma} d\Gamma
\end{equation}
Here,\( f(x) \) represents the image intensity function, and \( \mu(\Omega) \) is the mean of \( f(x) \) within the region \( \Omega \), defined by:
\begin{equation}\label{mu}
\begin{aligned}
\mu(\Omega) &= \frac{1}{|\Omega|} \int_{\Omega} f(x) \, dx, \quad &&\text{with} \quad |\Omega| = \int_{\Omega} dx, \\
\mu(\Omega^c) &= \frac{1}{|\Omega^c|} \int_{\Omega^c} f(x) \, dx.
\end{aligned}
\end{equation}
The first two terms \(E(\Omega)\) of the energy functional measure the variance of the image intensity within the segmented region \( \Omega \) and its complement \( \Omega^c \). The goal is to achieve a piecewise homogeneous segmentation, where each region exhibits minimal intensity variance around its respective mean value. Specifically, \( \mu(\Omega) \) and \( \mu(\Omega^c) \) represent the average intensities within the foreground and background regions, respectively. These terms enforce that the regions are homogeneous in intensity, which is desirable for accurate segmentation. 

The third term, introduces a regularization constraint on the boundary \( \Gamma \), controlling its length and smoothness. It helps to prevent over-segmentation and ensures that the boundary is well-defined, avoiding noisy or irregular boundaries. The parameter \( \eta \) is a weighting factor that balances the smoothness constraint against the variances, which represent data fidelity terms, allowing for flexibility in the segmentation process. The combined effect of these terms enables the active contour type models ~\cite{kashyap2017energy} to adapt the image structures while maintaining smooth boundaries and minimizing segmentation errors.

The main objective of the proposed model is to minimize the energy functional \( E(\Omega) \), which effectively aims to reduce the variance of the image intensity within both the segmented region \( \Omega \) and its complement \( \Omega^c \). By minimizing the first two terms of the energy functional, we achieve a segmentation where the regions are as homogeneous as possible in terms of their intensity. Specifically, the first term controls the intensity variance within \( \Omega \), while the second term ensures that the background region \( \Omega^c \) also exhibits minimal intensity variance. The mean intensity values \( \mu(\Omega) \) and \( \mu(\Omega^c) \), play a crucial role in this process. As the energy functional is minimized, the region \( \Omega \) evolves to better approximate homogeneous regions of the image, and similarly, the background \( \Omega^c \) becomes more homogeneous in terms of intensity.

Extending this framework to energy functionals, the shape derivative incorporates both domain and boundary integrals, reflecting geometric properties such as the domain's measure \( |\Omega| \), the measure of its complement \( |\Omega^c| \), and curvature effects represented by the mean curvature \( H \). The final expression for the shape derivative of the energy functional provides a comprehensive description of how domain variations influence energy distributions, accounting for both geometric and functional dependencies.

\section{Expression of the Shape Derivative}\label{Expression of the Shape Derivative}
Let \( \Omega \subset \mathbb{R}^N \) be a domain, where \( N \in \mathbb{N} \), and let \( T_t : \mathbb{R}^N \to \mathbb{R}^N \), with \( t \in [0, \varepsilon) \), represent a family of transformations that describes the deformation of \( \Omega \). The deformed domain is denoted by \( T_t(\Omega) = \Omega_t \).
The Eulerian velocity field ~\cite{delfour2011shapes} \( V(t, x) \)  at a point \( x(t) \) is defined as:
\begin{equation}
V(t, x) = \frac{\partial x}{\partial t} \left(t, T_t^{-1}(x)\right).
\end{equation}
The Eulerian derivative of a shape functional ~\cite{sokolowski1992introduction} \( J \) at \( \Omega \) in the direction of the vector field \( V \) is expressed as:
\begin{equation}
dJ(\Omega; V) = \lim_{t \downarrow 0} \frac{1}{t} \left( J(\Omega_t) - J(\Omega) \right).
\end{equation}
According to the structural theorem ~\cite{chicco2017shape}, if \( J \) is a shape functional that is differentiable at \( \Omega \) with respect to the velocity field \( V \), then there exists a functional \(\delta J  \in (C^k(\Gamma))' \) (referred to as the shape gradient) that satisfies:
\begin{equation}
dJ(\Omega; V) = \langle \delta J, v_n \rangle_{C^k(\Gamma)},
\end{equation}
\( v_n = V(0) \cdot n \). Furthermore, if the gradient, \( \delta J \in L^1(\Gamma) \), we have:
\begin{equation}
dJ(\Omega; V) = \int_{\Gamma} \delta J ~v_n \, d\Gamma.
\end{equation}

Consider the volume functional \( J(\Omega) = \int_{\Omega} dX \)~\cite{sokolowski1992introduction}, the shape derivative is given by
\begin{equation}
dJ(\Omega; V) = \int_{\Gamma} V(0) \cdot \hat{n} \, d\Gamma,
\end{equation}
\( \hat{n} \) is the outward unit normal vector to \( \Gamma\).\\
Let $Y:\mathbb{R}^N \to \mathbb{R}$, and shape functional defined by $J(\Omega)=\int_{\Omega}Y dx$ then shape derivative  is ~\cite{sokolowski1992introduction}:
\begin{equation}
dJ(\Omega; V) = \int_{\Gamma} Y \, V(0) \cdot \hat{n} \, d\Gamma.
\end{equation}
Suppose $f: \mathbb{R}^N \to \mathbb{R}^N$ and the functional defined as $J(\Omega)=\int_{\Gamma}f d\Gamma$, then shape derivative is defined by ~\cite{sokolowski1992introduction}
\begin{equation}
dJ(\Omega; V) = \int_{\Gamma} \left( \nabla f \cdot V(0) + f \, \text{div}_{\Gamma} V(0) \right) d\Gamma.
\end{equation}
\(\text{div}_{\Gamma} V(0)\) is the tangential divergence. Using this results, derived the shape derivative in next section, for the energy functional  $E(\Omega)$ in equation \ref{Mumford-shah model}.

\subsection{Shape Derivative of Mumford-Shah Energy}\label{Shape Derivative Mumford and Shah Model}
For simplicity, the proposed energy functional \( E(\Omega) \) is split into three components: \( E_{1}(\Omega) \), \( E_{1}(\Omega^c) \), and \( E_{2}(\Gamma) \).That is  
\begin{align}\label{Energy_combonents}
E_{1}(\Omega) &= \frac{1}{|\Omega|} \int_{\Omega} (f(x) - \mu(\Omega))^2 \, dx, \\
E_{1}(\Omega)^c &= \frac{1}{|\Omega^c|} \int_{\Omega^c} (f(x) - \mu(\Omega^c))^2 \, dx, \\
E_{3}(\Gamma)&= \int_{\Gamma} d\Gamma
\end{align}
Substituting equation \ref{mu} in \(E_{1}(\Omega)\) and \(E_{2} (\Omega)\) we obtain:
\begin{align}
E_{1}(\Omega) = \frac{1}{|\Omega|}\left( \int_{\Omega} f(x)^2 \, dx - \frac{1}{|\Omega|} \left( \int_{\Omega} f(x) \, dx \right)^2 \right) \label{Energy_combonents_with_expansion_1}\\
E_{1}(\Omega^c) = \frac{1}{|\Omega^c|}\left(  \int_{\Omega^c} f(x)^2 \, dx - \frac{1}{|\Omega^c|} \left( \int_{\Omega^c} f(x) \, dx \right)^2 \right) \label{Energy_combonents_with_expansion_2}
\end{align}
Applying the product and chain rule along with Gauss theorem, the shape gradient of  \(E_{1}(\Omega)\)  as follows: 
\begin{align}
\delta E_1(\Omega) &= \frac{1}{|\Omega|} f(x)^2 
- \frac{1}{|\Omega|^2} \int_\Omega f(x)^2 \, dx \notag \\
&\quad \quad \quad \quad \quad \quad \quad - \left\{ \frac{2f(x)}{|\Omega|^2} \int_\Omega f(x) \, dx 
- \frac{2}{|\Omega|^3} \left( \int_\Omega f(x) \, dx \right)^2 \right\} \notag \\
&= \frac{1}{|\Omega|} f(x)^2 
- \frac{1}{|\Omega|^2} \left( \int_\Omega f(x)^2 \, dx 
+ 2 f(x) \int_\Omega f(x) \, dx 
- \frac{2}{|\Omega|} \left( \int_\Omega f(x) \, dx \right)^2 \right)
\end{align}
Since \(\Omega^c=D-\Omega\), the equation \ref{Energy_combonents_with_expansion_2} will be,
\begin{align}
E_2(\Omega) &= \frac{1}{|D - \Omega|} \left\{ \int_D f(x)^2 \, dx - \int_{\Omega} f(x)^2 \, dx \right\} \quad \notag\\
&  \quad \quad \quad \quad \quad \quad \quad \quad - \frac{1}{|D - \Omega|^2} \left\{ \left( \int_D f(x)\,dx - \int_{\Omega} f(x)\,dx \right)^2 \right\} \notag \\
&= \frac{1}{|D - \Omega|} \int_D f(x)^2 \, dx - \frac{1}{|D - \Omega|} \int_{\Omega} f(x)^2 \, dx- \frac{1}{|D - \Omega|^2} \left( \int_D f(x)\,dx \right)^2 \notag \\
& \quad + \frac{2}{|D - \Omega|^2} \left( \int_D f(x)\,dx \int_{\Omega} f(x)\,dx \right) - \frac{1}{|D - \Omega|^2} \left( \int_{\Omega} f(x)\,dx \right)^2 
\end{align}
Shape gradient of \( E_{2}(\Omega)\) is given by:
\begin{align}
\delta E_{2}(\Omega) =\ & \frac{1}{|D - \Omega|^2} \int_D f(x)^2\ dx - \left( \frac{f(x)^2}{|D - \Omega|} + \frac{ \int_{\Omega} f(x)^2\,dx }{|D - \Omega|^2} \right) - \frac{2 \left( \int_D f(x)\,dx \right)^2}{|D - \Omega|^3} \notag  \\
&+ 2 \int_D f(x)\,dx  \left( \frac{f(x)}{|D - \Omega|^2}  + \frac{2}{|D - \Omega|^3} \int_D f(x)\,dx \right) \notag \\
& \quad \quad \quad \quad \quad \quad \quad - \left( \frac{2f(x)}{|D - \Omega|^2} \int_{\Omega} f(x)\,dx 
+ \frac{2}{|D - \Omega|^3} \left( \int_D f(x)\,dx \right)^2 \right)\notag  \\
=\ &\frac{1}{|\Omega^c|^2}\left( 
 \int_{\Omega^c} f(x)^2 \, dx 
- \frac{2}{|\Omega^c|} \left(
\left( \int_{D} f(x) \, dx \right)^2 
+ \left( \int_{\Omega} f(x) \, dx \right)^2+\int_{\Omega}f(x) 
\right)
\right)\notag \\
&\quad \quad \quad \quad \quad -\frac{1}{|\Omega^c|} f(x)^2\label{shape_gradient_E_2}
\end{align}
The third term \(E_{3}(\Gamma)\) specifically associated with the boundary of \( \Omega \) introduces a term involving the curvature \( H \), which represents the shape derivative of the boundary energy, so \(\delta E_{3}(\Gamma)=H\).
The shape gradient of proposed model is given by:
\begin{equation}
\delta E = \sum_{i=1}^{3} \delta E_i
\end{equation}
The velocity vector \( V(0) \) is typically chosen as the unit normal vector \(\hat{n}\) to the boundary \( \Gamma \) because normal displacements are the primary contributors to changes in the domain’s geometry. Additionally, using the normal direction naturally introduces geometric quantities such as the mean curvature \( H \), which plays a crucial role in energy functionals involving boundary regularization. From an optimization perspective, normal displacements correspond to the steepest descent direction, ensuring efficient energy minimization.

\section{Experimetal validation}\label{Experimental validation}
The algorithm is designed as an iterative process for boundary evolution. Initially, a polygonal boundary $\Gamma^{(0)}$ is initialized, represented by a set of vertices of $V_\Gamma=\{v_i:i=1...n\}$ with uniform edge lengths and a convergence threshold $E_{thr}$ is set to terminate the algorithm. The region enclosed by $\Gamma$ is denoted as $\Omega$, is computed as a binary mask where pixels inside the boundary are marked as foreground, and those outside as background. The mean intensities $\mu(\Omega)$ and $\mu(\Omega^c)$ were then calculated within the boundary ($\Omega$) and its complement ($\Omega^c$), respectively. Optionally, the energy functional $E(\Omega)$ could be computed to verify if the segmentation quality had improved as the method progressed.

Subsequently, the gradient \( \delta E \) is computed at each vertex of \( \Gamma \) to determine the direction vector for each point of $\Gamma$. The boundary is then updated by moving each vertex in the direction of the negative normal of length $\delta E$, with a step size \( dt \).
The boundary is updated to $(i+1)^{th}$ state from its $i^{th}$ state as follows :
\begin{equation}
\Gamma^{(i+1)} = \Gamma^{(i)} - dt (\delta E)^{(i)}
\end{equation}
To maintain a consistent resolution, the boundary is resampled periodically, redistributing the vertices evenly along the curve. Convergence is checked by comparing to the predefined threshold $E_{thr}$; if relative difference of energy fell below this threshold, the algorithm is terminated, assuming that further updates would not significantly improve the segmentation. Otherwise, the process is repeated from the computation of the binary mask until convergence is achieved.

\section{Results analysis}\label{Result analysis}
\begin{figure}[H]
    \centering
    \begin{minipage}{0.22\linewidth}
        \centering
        \includegraphics[width=\linewidth]{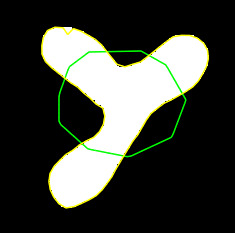}
    \end{minipage}
    \begin{minipage}{0.22\linewidth}
        \centering
        \includegraphics[width=\linewidth]{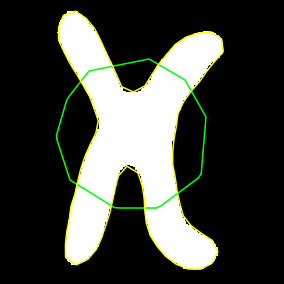}
    \end{minipage} 
    \begin{minipage}{0.23\linewidth}
        \centering
        \includegraphics[width=\linewidth]{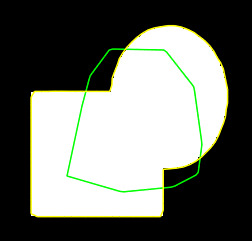}
    \end{minipage}
    \begin{minipage}{0.215\linewidth}
        \centering
        \includegraphics[width=\linewidth]{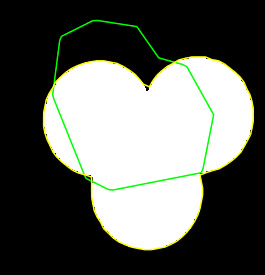}
    \end{minipage} \\[0.5em]

\begin{minipage}{0.22\linewidth}
        \centering
        \includegraphics[width=\linewidth]{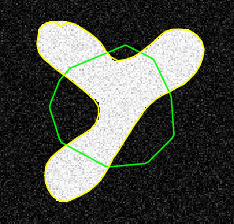} \caption*{(a)}
    \end{minipage}
    \begin{minipage}{0.22\linewidth}
        \centering
        \includegraphics[width=\linewidth]{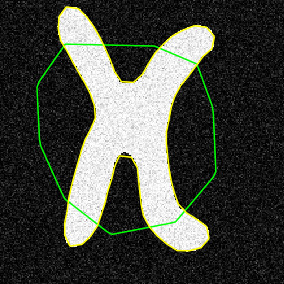} \caption*{(b)}
    \end{minipage} 
    \begin{minipage}{0.215\linewidth}
        \centering
        \includegraphics[width=\linewidth]{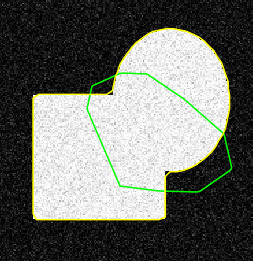} \caption*{(c)}
    \end{minipage}
    \begin{minipage}{0.22\linewidth}
        \centering
        \includegraphics[width=\linewidth]{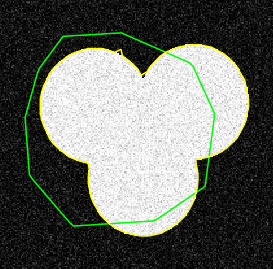}\captionsetup{labelformat=empty} 
        \caption*{(d)}
      \end{minipage}
\caption{Segmentation of synthetic images. The first row shows binary images, and the second row shows images with Gaussian noise ($SD=25$). The green polygon represents the initial boundary with $100$ points, and the yellow polygon represents the final boundary after $100$ iterations}
\label{img_synthetic}
\end{figure}

We conducted experiments on a variety of synthetic images, including binary images, images corrupted with Gaussian noise (see Figure \ref{img_synthetic}), and color images in both RGB and LAB spaces. The segmentation process begins with an initial boundary set of vertices, depicted in green color, and iteratively evolves it to the final boundary as shown in yellow color. The energy functional guiding this evolution consists of three key terms: (1) Variance of pixel intensities inside the boundary, (2) Variance of pixel intensities outside the boundary, and (3) A boundary regularization term. The first term ensures homogeneity within the segmented region, the second term ensures differentiation from the background, and the third term helps maintain perimeter and smoothness in the boundary evolution.

\begin{figure}[H]
    \centering
    \begin{minipage}{0.45\linewidth}
        \centering
        \includegraphics[width=\linewidth]{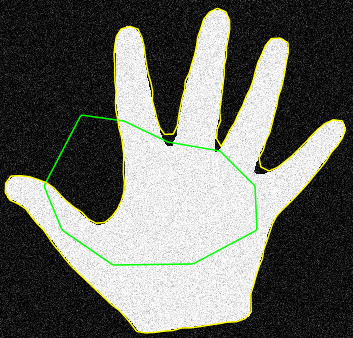}\caption*{(a)}
    \end{minipage}
    \begin{minipage}{0.40\linewidth}
        \centering
        \includegraphics[width=\linewidth]{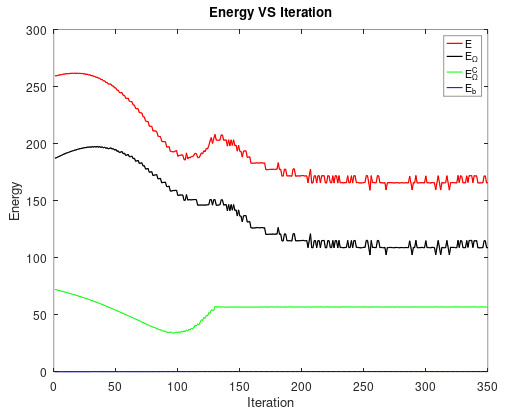}  \caption*{(b)}
    \end{minipage}
   \caption{(a) Segmentation of a noisy palm image with a resolution of $250 \times 250$ pixels. The green polygon represents the initial  boundary with 150 points, while the yellow polygon shows the converged result after 350 iterations. (b) Energy curve illustrating the different energy terms of the functional. The red, black, green, and blue curves represent the total energy, $E_{\Omega}$, $E_{\Omega^c}$, and $E_{{\text{boundary}}}$, respectively.}
   \label{im_palm}
\end{figure}

It is observed that the variance of pixel intensities inside and outside the initial boundary progressively decreased as the segmentation process advanced. This reduction in variance indicates that the boundary effectively separates regions of distinct intensities, ultimately converging towards the true boundary of the object. The performance of the method is evaluated using various types of images, including clean and noisy images, demonstrating its robustness against noise and its applicability to grayscale images.
In Figure~\ref{im_palm}(a), an image of a human palm with added Gaussian noise (standard deviation = 25) is analyzed for boundary evolution. The initial boundary, marked in green, is set at the beginning of the segmentation process, whereas the evolved boundary, depicted in yellow, is obtained after 350 iterations. This transformation illustrates how the boundary effectively adjusted to the palm’s shape despite the noise interference. In Figure~\ref{im_palm}(b), the energy evolution over iterations is plotted, where different energy components are represented using distinct colors. The total energy (\(E\)) is shown in red, while the first, second, and third energy terms (\(E_{1}\), \(E_{2}\), and \(E_{3}\)) are represented in black, green, and blue, respectively. Over the course of the iterations, the total energy is observed to decrease and stabilize after approximately 200 iterations, indicating the convergence of the segmentation process.

Furthermore, the segmentation method successfully captured the boundaries of various objects, even in the presence of significant noise. The final boundary accurately aligned with object edges, validating the precision of the proposed approach. This confirms that minimizing the variance within and outside the boundary is a powerful strategy for effective image segmentation.

\begin{figure}[H]
    \centering
    \begin{minipage}{0.35\linewidth}
        \centering
        \includegraphics[width=\linewidth]{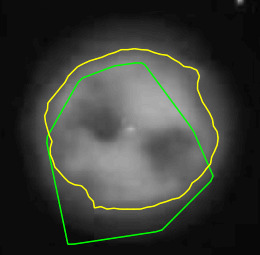}\caption*{(a)}
    \end{minipage}
    \begin{minipage}{0.50\linewidth}
        \centering
        \includegraphics[width=\linewidth]{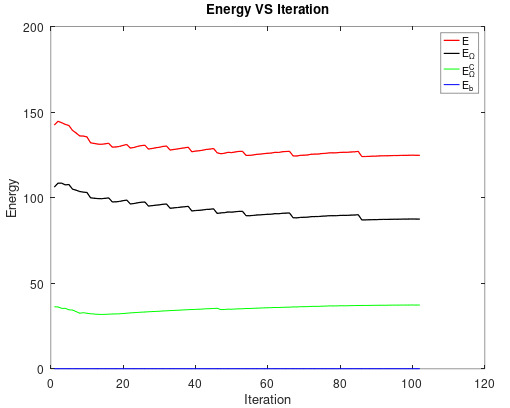}\caption*{(b)}
    \end{minipage} 
\caption{(a) Segmentation of a galaxy image (250 × 250 pixels) with initial (green) and final (yellow) polygons after 100 iterations. (b) Energy plot showing the total energy (\(E\), red) and its components: \(E_{1}\) (black), \(E_{2}\) (green), and \(E_{3}\) (blue)}  \label{fig_galaxy}
\end{figure}

In Figure~\ref{fig_galaxy} (a), a grayscale image of a galaxy is used. Since the galaxy does not have a well-defined boundary, the initial polygon (marked in green) is placed arbitrarily around the central bright region. Over the course of 100 iterations, the boundary evolved and adapted to enclose the galaxy's core, as shown by the yellow polygon. The energy evolution plot in Figure~\ref{fig_galaxy} (b) shows that the total energy (\(E\), red) and its individual components gradually decrease and stabilize, indicating the convergence of the segmentation process. This result demonstrates that the model successfully detects the main structure of the galaxy despite the absence of a clear boundary.  

\begin{figure}[H]
    \centering
    \begin{minipage}{0.45\linewidth}
        \centering
        \includegraphics[width=\linewidth]{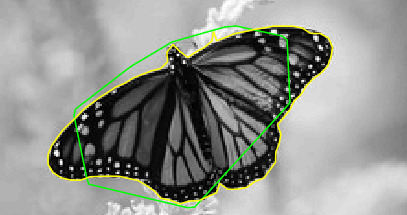}\caption*{(a)}
    \end{minipage}
    \begin{minipage}{0.40\linewidth}
        \centering
        \includegraphics[width=\linewidth]{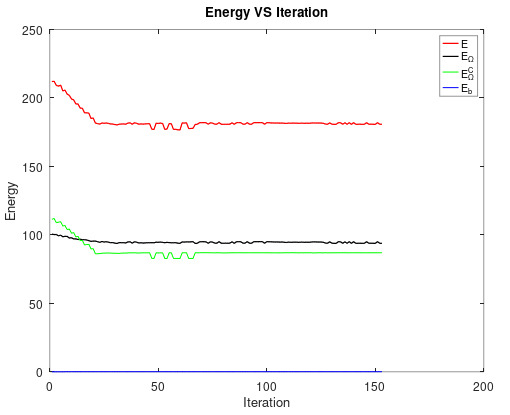}\caption*{(b)}
    \end{minipage}   
    \caption{(a) Segmentation of a butterfly image (248 × 150 pixels) with initial (green) and final (yellow) polygons after 100 iterations. (b) Energy plot showing the total energy (\(E\), red) and its components: \(E_{1}\) (black), \(E_{2}\) (green), and \(E_{3}\) (blue)}
    \label{fig:butterfly}
\end{figure}

Figure~\ref{fig:butterfly} (a) presents a grayscale image of a butterfly. The initial boundary (green) is placed around the object, and after 25 iterations, it converged to the butterfly’s natural boundary, as shown in yellow. The energy plot in Figure~\ref{fig:butterfly} (b) shows a steady decrease in total energy, stabilizing after approximately 100 iterations. Unlike the galaxy image, the butterfly has a well-defined shape, which allows the polygon to precisely capture its edges. The stabilization of energy terms confirms that the segmentation process effectively extracts the object's boundaries. These results illustrate the adaptability of the method to different types of grayscale images, whether they have clear or ambiguous boundaries.

\subsection{Color Segmentation in RGB Space}
\begin{figure}[H]
    \centering
    \begin{minipage}{0.45\linewidth}
        \centering
        \includegraphics[width=\linewidth]{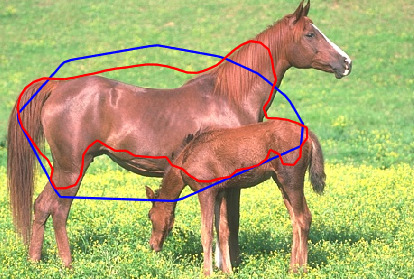}
    \end{minipage}
    \begin{minipage}{0.45\linewidth}
        \centering
        \includegraphics[width=\linewidth]{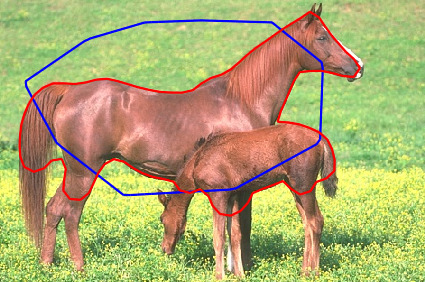}
    \end{minipage}\\
    
    \centering
    \begin{minipage}{0.47\linewidth}
        \centering
        \includegraphics[width=\linewidth]{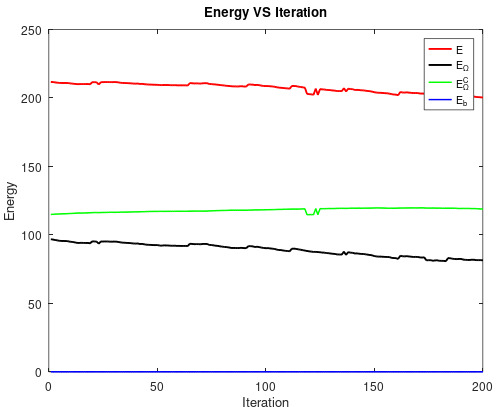}\caption*{(a)}
    \end{minipage}
    \begin{minipage}{0.47\linewidth}
        \centering
        \includegraphics[width=\linewidth]{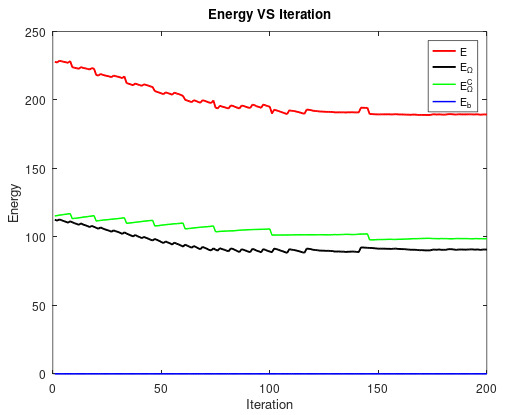}\caption*{(b)}
    \end{minipage}

    \caption{(a) Segmentation in RGB space. (b) Segmentation in \textit{L\textsuperscript{*}a\textsuperscript{*}b\textsuperscript{*}} space initial curve (blue) and final (red) curve  after 200 iterations. }
    \label{fig:Hourse}
\end{figure}

Image segmentation was performed directly in the RGB color space to fully utilize the raw color information available in standard image formats (see Figure \ref{fig:Hourse} (a)). The original image was first decomposed into its red, green, and blue channels. For each channel, the shape derivative of the energy functional was computed independently. The resulting channel-wise gradients were then summed to form a single descent direction, which guided the evolution of the contour via gradient descent. The initial set of vertices was placed around the region of interest, represented by a closed polygonal curve in blue. After 200 iterations, the set of points evolved to align closely with the object boundary, forming the red curve. We also conducted experiments in the LAB color space, which provided even better segmentation results due to its perceptual uniformity.
\subsection{Color Segmentation in \textit{L\textsuperscript{*}a\textsuperscript{*}b\textsuperscript{*}} Space}
\begin{figure}[H]
    \centering
    \begin{minipage}{0.30\linewidth}
        \centering
        \includegraphics[width=\linewidth]{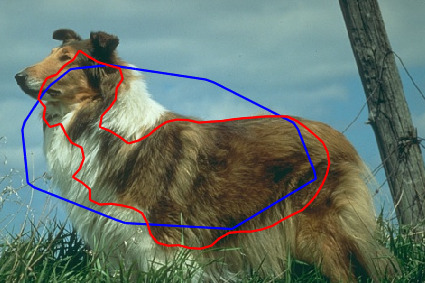} \\
       
    \end{minipage}
    \begin{minipage}{0.30\linewidth}
        \centering
        \includegraphics[width=\linewidth]{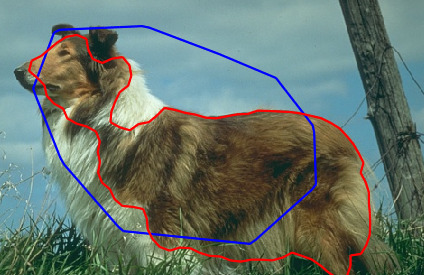} \\
      
    \end{minipage}   
    \begin{minipage}{0.30\linewidth}
        \centering
        \includegraphics[width=\linewidth]{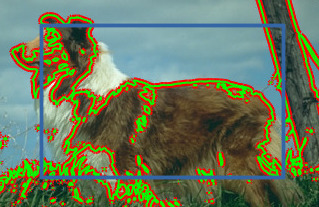} \\
       
    \end{minipage}

    \begin{minipage}{0.30\linewidth}
        \centering
        \includegraphics[width=\linewidth]{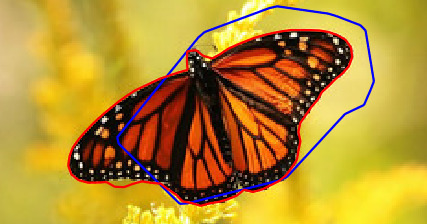} \\
        (a)
    \end{minipage}
    \begin{minipage}{0.30\linewidth}
        \centering
        \includegraphics[width=\linewidth]{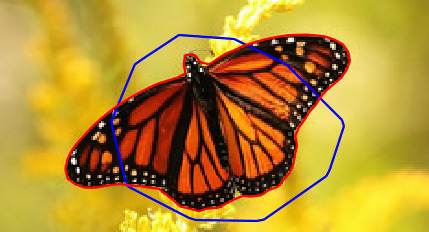} \\
        (b)
    \end{minipage}   
    \begin{minipage}{0.30\linewidth}
        \centering
        \includegraphics[width=\linewidth]{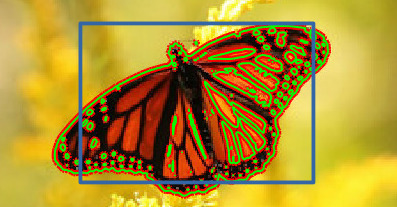} \\
        (c)
    \end{minipage}
    
    \caption{(a) Segmentation in RGB space (b) Segmentation in LAB space using proposed model. (c) Segmentation in RGB space using the Chan-Vese model. initial curve (blue) and final (red) curve.}
    \label{fig:RGB_LAB_CV}
\end{figure}

The image segmentation was performed using the LAB color space to take advantage of its perceptual uniformity and better color separation. The original RGB image was converted into LAB, and segmentation was carried out by evolving an initial contour based on the gradient descent of a shape energy functional. The shape derivatives were computed separately for each of the L, A, and B channels, and their combination guided the contour evolution. Over several iterations, the contour adapted to the object boundaries in the image, with periodic resampling to maintain stability. This approach allowed for more accurate segmentation of regions with subtle color differences, which are often hard to distinguish in the RGB space.
In the RGB space (Figure (a)), the segmentation roughly captures the shapes of the horses but fails to accurately align with finer details, especially around complex boundaries. In contrast, when using the LAB color space (Figure (b)), the segmentation significantly improves, producing contours that closely follow the true object boundaries. This improvement can be attributed to the perceptual uniformity of the LAB space, which better separates color information relevant for distinguishing different regions. While the RGB-based implementation benefits from its direct use of standard image data, it proves sensitive to illumination changes and often yields irregular boundaries in areas of subtle color variation. On the other hand, the LAB-based approach, by decoupling luminance from chromatic components and leveraging a perceptually uniform distance metric, consistently produces cleaner, more stable, and more accurate segmentations. These results demonstrate that applying the model in the LAB space leads to superior segmentation performance compared to the RGB space. The figure \ref{fig:RGB_LAB_CV} show the segmentation of two image with in RGB channel and LAB channel with proposed model. colums (c) repressent the segemtation of smane images with chanvese model. due to level set frame work the intitial curve is  splinted into  multiple.

Figure \ref{fig:RGB_LAB_CV} shows the segmentation of two images in both the RGB and LAB color spaces using the proposed model. Columns (c) represent the segmentation of same images using the Chan-Vese model. Due to the level set framework, the initial curve is split into multiple.
The observed convergence behavior and energy profiles demonstrate the effectiveness of the proposed framework, making it a promising technique for various image processing applications. This approach effectively segments images by employing the Mumford-Shah energy model without relying on the level sets or parametrized curves. The results indicate that the method is highly efficient, robust to noise, and capable of accurately delineating object boundaries in grayscale images.

\section{Conclusion}
We proposed a modified non-parametric level set free, Mumford-Shah model for image segmentation using shape gradients. This method effectively minimizes texture-type intensity variance both inside and outside the polygon of boundary points, ensuring precise object boundary detection, as demonstrated by experimental analysis. Extensive experiments on synthetic and noisy images highlighted its robustness and efficiency, with the boundary evolving smoothly to align with object contours. Additionally, we demonstrated the method's effectiveness on color images in both RGB and LAB channels. Energy analysis confirmed the method's efficacy, showing consistent variance reduction and faster convergence compared to traditional techniques. The approach proved to be computationally efficient, resilient to noise, and successful in segmenting both grayscale and color images, making it well-suited for real-world applications.

\subsection*{Data Availability}
  The authors affirm that the data used for experimenting the model will be made available if needed.

\subsection*{Grammar and Readability Disclosure}
  This document has been reviewed with AI-based tools \cite{chatgpt} to check grammar and address readability improvements.

\bibliographystyle{plain}
\bibliography{references}

\end{document}